\documentclass[11pt]{article}
\usepackage{epsfig}
\usepackage{color,epsfig}
\usepackage{amsfonts,graphicx,psfrag}
\usepackage{float}
\usepackage{graphicx}

\def\be{\begin{equation}}
\def\ee{\end{equation}}
\def\bea{\begin{eqnarray}}
\def\eea{\end{eqnarray}}
\def\bes{\begin{eqnarray*}}
\def\ees{\end{eqnarray*}}

\def\nn{\nonumber}
\def\<{\langle}
\def\>{\rangle}
\def\lb{\label}
\def\bs{\setminus}

\def\R{{\bf R}}
\def\C{{\bf C}}
\def\Z{{\bf Z}}
\def\N{{\bf N}}
\def\U{{\bf U}}
\def\Q{{\bf Q}}
\def\T{{\bf T}}

\def\ga{{\gamma}}

\def\th{{\theta}}
\def\om{{\omega}}
\def\Om{{\Omega}}
\def\ep{{\epsilon}}
\def\lm{{\lambda}}

\def\sg{{\sigma}}

\def\Sg{{\Sigma}}

\def\H{{\cal H}}
\def\T{{\cal T}}

\def\P{{\cal P}}
\def\J{{\cal J}}

\def\Sp{{\rm Sp}}

\def\dm{{\rm \diamond}}

\def\ol#1{\overline{#1}}  

\def\ol#1{\overline{#1}}  

\title{Stability of closed characteristics on compact convex \\
hypersurfaces in $\R^{2n}$}
\author{Xijun Hu\thanks{Partially supported
by NSFC(No.11131004), PCSIRT ( IRT1264) and NCET, E-mail:xjhu@sdu.edu.cn }
\quad Yuwei Ou \thanks{Partially supported
by NSFC(No.11131004), E-mail:yuweiou@163.com }\\ \\
Department of Mathematics, Shandong University\\
Jinan, Shandong 250100, The People's Republic of China\\
}
\date{}

\date{}

\begin{document}

\maketitle

\begin{abstract}
{Let $\Sigma\subset \R^{2n}$ with $n\geq2$ be any $C^2$ compact convex hypersurface and  only has  finitely geometrically distinct closed characteristics.   Based on Y.Long and C.Zhu 's  index jump  methods \cite{LoZ1}, 
we prove that  there are at least two  geometrically distinct elliptic closed characteristics,  and   moreover,   there exist at least $\varrho_{n} (\Sigma)$ ($\varrho_{n}(\Sigma)\geq[\frac{n}{2}]+1$)
 geometrically distinct closed characteristics
such that for any two elements among them,
the ratio of their mean indices is irrational number.  }
\end{abstract}



{\bf Running head}: Stability of closed characteristics

\renewcommand{\theequation}{\thesection.\arabic{equation}}
\renewcommand{\thefigure}{\thesection.\arabic{figure}}

\setcounter{equation}{0}

\section{Introduction and main results}
In this paper, let $\Sigma$ be a fixed $C^2$ compact convex hypersurface
in $\R^{2n}$, i.e., $\Sigma$ is the boundary of a compact and strictly
convex region $U$ in $\R^{2n}$. We denote the set of all such hypersurfaces
by $\H(2n)$. Without loss of generality, we suppose $U$ contains the origin.
We consider closed characteristics $(\tau, x)$ on $\Sigma$, which are
solutions of the following problem
\bea
\cases{\dot{x}(t) &$= JN_{\Sigma}(x(t)),
\quad x(t)\in\Sigma, \qquad \forall t\in\R,$ \cr
x(\tau) &$= x(0),$ \cr}
\label{6.20}\eea
where $J=\left(\matrix{0 &-I_n\cr
                I_n  & 0\cr}\right)$,
$I_n$ is the identity matrix in $\R^n$, $\tau>0$, $N_\Sigma(x)$ is
the outward normal vector of $\Sigma$ at $x$ normalized by the
condition $N_\Sigma(x)\cdot x=1$. Here $a\cdot b$ denotes the
standard inner product of $a, b\in\R^{2n}$. A closed characteristic
$(\tau, x)$ is {\it prime}, if $\tau$ is the minimal period of $x$.
Two closed characteristics $(\tau, x)$ and $(\sigma, y)$ are {\it
geometrically distinct},  if $x(\R)\not= y(\R)$. We denote by $\T(\Sg)$, the
set of all geometrically distinct closed characteristics on $\Sigma$,
and $[(\tau,x)]$ the set of all closed characteristics which are geometrically
the same as $(\tau,x)$. $^{\#}$A denotes the total number of elements in a set A.

Let $j:\R^{2n}\rightarrow \R$ be the gauge function of $\Sigma$, i.e.,
$j(\lambda x)=\lambda$ for $x\in\Sigma$ and $\lambda\geq0$, then
$j\in C^{2}(\R^{2n}\setminus \{0\},\R)\cap C^{1}(\R^{2n},\R)$ and  $\Sigma=j^{-1}(1)$. Fix a constant
$\alpha\in(1,2)$ and define the Hamiltonian function $H_{\alpha}:\R^{2n}\rightarrow [0,+\infty]$
by
\be
H_{\alpha}(x)=j(x)^{\alpha},  \ \ \ \ \forall x\in \R^{2n}. \label{1.2}
\ee
Then $H_{\alpha}\in C^{2}(\R^{2n}\setminus\{0\},\R)\cap C^{1}(\R^{2n},\R)$
is convex and $\Sigma=H^{-1}_{\alpha}(1)$.
It is well known that the problem (1.1) is equivalent to the following given energy
problem of the Hamiltonian system
\bea
\cases{\dot{x}(t) &$= JH_{\alpha}^{\prime}(x(t)),
\quad H_{\alpha}(x(t)) =1, \qquad \forall t\in\R,$ \cr
x(\tau) &$= x(0).$ \cr}
\label{1.3}\eea
Denote by $\J(\Sigma,\alpha)$ the set of all geometrically distinct solutions $(\tau,x)$ of (1.3) where $\tau$
is the minimal period of $x$. Note that elements in
$\T(\Sigma)$ and $\J(\Sigma,\alpha)$ are
one to one correspondent to each other.

Let $(\tau,x)\in \J(\Sigma,\alpha)$, the fundamental solution $\gamma_{x}: [0,\tau]\rightarrow \Sp(2n)$ with
$\gamma_{x}(0)=I_{2n}$ of the linearized Hamiltonian system
\be
\dot{y}(t)=JH^{\prime\prime}_{\alpha}(x(t))y(t), \ \ {\rm for\;\ all}\;\ t\in \R \\
\ee
is called the associated symplectic path of $(\tau,x)$. The eigenvalues of $\gamma_{x}(\tau)$ are called
Floquet multipliers of $(\tau,x)$.
A closed characteristic $(\tau,x)$ is {\it non-degenerate}, if $1$ is a Floquet
multiplier of $y$ of precisely algebraic multiplicity $2$, and is {\it elliptic},
if all the Floquet multipliers of $x$ are on ${\bf U}=\{z\in\C\,|\,|z|=1\}$,
i.e., the unit circle in the complex plane.

The study on closed characteristics on the star-shaped hypersurface in the global sense started in 1978 by  Rabinowitz in \cite{Rab1} and Weinstein for the convex hypersurface independently\cite{Wei1}.
For more results on the multiplicity of geometrically distinct closed characteristics on convex hypersurfaces, please refer to
\cite{EkL1}, \cite{EkH1}, \cite{Szu1}, \cite{HWZ1}, \cite{LoZ1}, \cite{LLZ1}, \cite{WHL}, etc. and references therein.

For the stability problem,  Ekeland proved in [Eke2] the existence of at least one elliptic closed
characteristic on $\Sigma$ provided $\Sigma \in \H(2n)$ is $\sqrt{2}$-pinched. In \cite{DDE1} of 1992, Dell'Antonio, D'Onofrio
and Ekeland proved the existence of at least one elliptic closed characteristic on $\Sigma$ provided $\Sigma \in \H(2n)$
satisfies $\Sigma=-\Sigma$. In \cite{Lon4} of 2000, Long proved that $\Sigma \in \H(4)$ and $^{\#}\T(\Sigma)=2$ imply that
both of the closed characteristics must be elliptic. In \cite{LoZ1} of 2002, Long and Zhu further proved
when $^{\#}\T(\Sigma)<+\infty$, there exists at least one elliptic closed characteristic and there are at least
$[\frac{n}{2}]$ geometrically distinct closed characteristics on $\Sigma$ possessing irrational mean indices, which are
then non-hyperbolic. Moreover,  they proved there exist at least two elliptic closed characteristics provided that $^{\#}\T(\Sigma)\leq2\varrho_n(\Sigma)-2<\infty$ ,  where $\varrho_n(\Sigma)$ is defined by Definition 2.10.  In the recent paper \cite{LoW1}, Long and Wang proved that there exist at least
two non-hyperbolic closed characteristic on $\Sigma \in \H(6)$ when $^{\#}\T(\Sigma)<+\infty$ and in \cite{Wa1}, Wang proved that
there exist at least two elliptic closed characteristic on $\Sigma \in \H(6)$ when $^{\#}\T(\Sigma)=3$.   Other results please refer \cite{Lon2}, \cite{Lon4}, \cite{LLW}.  Motivated by these
results, we prove the following results in this paper:
\\\\
{\bf Theorem 1.1.} For any
$\Sigma\in\H(2n)$ with $n\geq2$ satisfying $^{\#}\T(\Sigma)<+\infty$,
there exist at least two elliptic closed characteristics on $\Sigma$.
\\\\
A typical example is the non-resonant  ellipsoid in $\R^{2n}$, that is $\Sigma$ is defined by
\bea \sum_{i=1}^n \frac{\alpha_i}{2} (p_i^2+q_i^2)=1, \eea
where $\alpha_i/\alpha_j\in\R\setminus  \Q$.   There just exist $n$  closed characteristics $x_i,i=1,...,n$ and  their mean Maslov-type index  satisfy $\hat{i}(x_i)/\hat{i}(x_j)=\alpha_j/\alpha_i\in\R\setminus\Q$. When$^{\#}\T(\Sigma)$  is finite, it seems that all the  Maslov-type index of the closed characteristics are similar to those in the non-resonant ellipsoid. Another example is  in the case $n=2$,
it has been proved in \cite{HWZ1} that there are either infinite or $2$ closed characteristics. When $n=2$, in the case $^{\#}\T(\Sigma)=2$,
 Long, Wang and Hu have proved that both of their mean index are irrational number and all the iterations of their Maslov-type index are the same as a non-resonant  ellipsoid, another different proof is given by \cite{BCE} . As results toward this aspect, we proved that
\\\\
{\bf Theorem 1.2.} For any
$\Sigma\in\H(2n)$ satisfying $^{\#}\T(\Sigma)<+\infty$,
there exist at least $\varrho_{n}(\Sigma)$
geometrically distinct closed characteristics on $\Sigma$
such that any two element $[(\tau,x)]$, $[(\tilde{\tau},\tilde{x})]$
satisfy
\bea
\frac{\hat{i}(x,1)}{\hat{i}(\tilde{x},1)}\in \R\setminus\Q.
\eea
\\
\\
The main ingredient in our proof of this theorems is
the Maslov-type index iteration theory developed by Long and his
coworkers, especially based on some new observations on the common index jump theorem of
Long and Zhu (Theorem 4.3 of \cite{LoZ1}, cf. Theorem 11.2.1
of \cite{Lon5}).  In Section 2,  we review briefly the variational
structure and the common index jump theorem of Long and Zhu with some further discussion, we prove our main theorems in Section 3. For reader's convenience, we brief review the Maslov-type index theory in Section 4.

In this paper, let $\N$, $\Z$, $\Q$, $\Q^{+}$, $\R$, and $\R^+$ denote
the sets of natural integers, integers, rational
numbers, positive rational number, real numbers, and positive real numbers respectively.
Denote by  $(a,b)$  and $|a|$ the standard inner product and norm in
$\R^{2n}$. Denote by $\langle\cdot,\cdot\rangle$ and $\|\cdot\|$
the standard $L^2$ inner product and $L^2$ norm. We also define the functions
\bea
\cases{[x]&$=\max\{k\in\Z|k\leq x\},
\quad E(x)=\min\{k\in\Z|k\geq x\},$  \cr
\{x\}&$=x-[x],\quad\varphi(x)=E(x)-[x].$ \cr} \label{1.2}
\eea

\setcounter{equation}{0}
\section{Brief review of Long-Zhu index jump Theorem and with further discussion}
To solve the given fixed energy problem (\ref{1.3}) as
in \cite{Eke3} instead, we consider the following
fixed period problem:
\bea
\cases{\dot{z}(t) &$= JH_{\alpha}^{\prime}(z(t)),
 \qquad \forall t\in\R,$ \cr
z(1) &$= z(0).$ \cr}
\label{6.21}\eea
Define
\be E=
\left\{u\in L^{(\alpha-1)/\alpha}(\R/\Z,\R^{2n})\left|\frac{}{}\right.\int_0^1u(t)dt=0\right\}.
 \label{2.9}\ee
The corresponding Clarke-Ekeland dual action function $f:E\rightarrow \R$ is defined by
\be f(u)=\int_0^1\left(\frac{1}{2}(Ju, Mu)+H^{*}_\alpha(-Ju)\right)dt,
\lb{2.11}\ee
where $Mu$ is defined by $\frac{d}{dt}Mu(t)=u(t)$ and $\int_0^1Mu(t)dt=0$, and the usual dual
function $H^{*}_{\alpha}$ of $H_{\alpha}$ is defined by
\bea
H^{*}_{\alpha}(x)=\sup_{y\in\R^{2n}}((x,y)-H_{\alpha}(y)),
\eea
 then $f\in C^{2}(E,\R)$.
Suppose $z$
is a solution of (\ref{6.21}), then $u=\dot{z}$ is a critical
point of $f$. Conversely, suppose $u\in E\setminus\{0\}$ is a critical point
of $f$, then there exists a unique $\xi_{u}\in\R^{2n}$ such that
$z_{u}(t)=Mu(t)+\xi_{u}$ is a solution of (\ref{6.21}). In particular, solutions
of (\ref{6.21}) are in one to one correspondence with critical
points of $f$.

Following $\S V.3$ of \cite{Eke3}, we denote by "ind" the Fadell-Rabinowitz $S^{1}$-
action cohomology index theory (please refer \cite{FaR1}) for $S^{1}$-invariant subsets of $E$ defined in \cite{Eke3}.
For $[f]_{c}\equiv\{u\in E|f(u)\leq c\}$, the following critical values of $f$ are defined
\bea
c_{k}=\inf\{c<0|{\rm ind}([f]_{c})\geq k\}, \ \ \forall k\in \N
\eea
Ekeland and Hofer proved the following theorem which is a basis of further study \cite{EkH1}, \cite{Eke3},  the theorem with the following form is from  \cite{LoZ1}.
\\
\\
{\bf Theorem 2.1.}
\bea  \nonumber
-\infty<c_{1}&=& \inf_{u\in E}f(u)\leq c_{2}\leq \cdots\leq c_{k}\leq c_{k+1}\leq\cdots <0 \\ \nonumber
c_{k}&\rightarrow& 0 \ \ {\rm as} \ \ k\rightarrow +\infty, \\ \nonumber
^{\#}\T(\Sigma)&=&+\infty \ \ {\rm if} \ \ c_{k}=c_{k+1} \ \ {\rm for \ \ some} \ \ k\in\N.
\eea
For any given $k\in\N$, there exists $(\tau,x)\in\J(\Sigma,\alpha)$ and $m\in \N$ such that for
\bea
u^{x}_{m}(t)&=&(m\tau)^{(\alpha-1)/(2-\alpha)}\dot{x}(m\tau t), \ \ 0\leq t\leq1, \label{2.6}
\eea there hold
\bea
f^{\prime}(u_{m}^{x})&=&0,\ \ f(u^{x}_{m})=c_{k}, \label{2.7}\\
i(x,m)&\leq&2k-2+n\leq i(x,m)+\nu(x,m)-1, \label{2.8}
\eea
where $i(x,m)$ is the Maslov-type index of closed characteristics $x$ with $m$-th iteration,  for the readers whom are not familiar with the Maslov-type index, please refer Section 4 for the definitions and basic notations.
\\\\
{\bf Definition 2.2.}(cf.\cite{LoZ1}) For any $\Sigma\in\H(2n)$ and $\alpha\in(1,2)$, $(\tau,x)\in \J(\Sigma,\alpha)$ is
$(m,k)$-variationally visible, if there exist some $m$ and $k\in \N$ such that (\ref{2.7}), (\ref{2.8}) hold for
$u_{m}^{x}$ defined by (\ref{2.6}). We call $(\tau,x)\in\J(\Sigma,\alpha)$ infinite variationally visible,
if there exist infinitely many $(m,k)$ such that $(\tau,x)$ is $(m,k)$-variationally visible. We denote
by $\mathcal{V}_{\infty}(\Sigma,\alpha)$ the subset of $\J(\Sigma,\alpha)$ in which
a representative $(\tau,x)\in\J(\Sigma,\alpha)$ of each $[(\tau,x)]$ is infinitely  variationally visible.
\\\\
{\bf Theorem 2.3.} (cf.\cite{LoZ1})  There  exists an integer $K\geq 0$ and an injection map
$$p : \N + K\longrightarrow \mathcal{V}_{\infty}(\Sigma, \alpha) \times \N$$ such that\\
(i) For any $k\in\N+ K$, $(\tau,x) \in \J(\Sigma,\alpha)$ and $m \in \N$ satisfying $p(k) = ([(\tau, x)], m)$, (\ref{2.7}) and (\ref{2.8}) hold, and
\\(ii) For any $k_j\in \N + K$, $k_1 < k_2$, $(\tau_j,x_j) \in \J(\Sigma,\alpha)$ satisfying $p(k_j) = ([(\tau_j, x_j)], m_j)$  with $j = 1, 2$,
\bea \hat{i}(x_1, m_1) < \hat{i}(x_2, m_2).\lb{2.9}   \eea
\\\\
{\bf Remark 2.4}
Since $\alpha\in(1,2)$,  Let $\ga\in\P_{\tau}(2n)$ be the fundamental matrix of $(x,\tau)$, $M=\ga(\tau)$,  then  there exist $P\in\Sp(2n)$ and $Q\in\Sp(2n-2)$ such that
$M=P^{-1}(N_{1}(1,1)\diamond Q)P$.   Since $H_\alpha$  is convex, then $i_1(\ga)\geq n$,  and consequently the   mean Maslov-type  index  $\hat{i}(\ga)>2$ for $n\geq2$.
\\\\
A key Theorem of \cite{LoZ1} is the following index jump theorem.
\\\\
{\bf Theorem 2.5.}(cf. P.350 of \cite{LoZ1}) Let $\ga_{k}\in\P_{\tau_{k}}(2n)$
for $k=1,\ldots q$ be a finite collection of symplectic paths. Let $M_{k}=\ga_{k}(\tau_{k})$.
Suppose that there exists $P_{k}\in\Sp(2n)$ and $Q_{k}\in\Sp(2n-2)$ such that
$M_{k}=P^{-1}_{k}(N_{1}(1,1)\diamond Q_{k})P_{k}$ and $\hat{i}(\ga_{k},1)>0$, for all $k=1,\ldots,q$
. Then there exist infinitely many $(N,m_{1},\ldots,m_{q})\in\N^{q+1}$ such that
\bea
I(k,m_{k})&=&N+\Delta_{k},\label{3.30}
\eea
where
\bea \nonumber
I(k,m_{k})&=&m_{k}(i(\gamma_{k},1)+S_{M_{k}}^{+}(1)-C(M_{k}))\\
&+&\sum_{\theta\in(0,2\pi)}E(\frac{m_{k}\theta}{\pi})
S_{M_{k}}^{-}(e^{\sqrt{-1}\theta}),\label{3.31}\\
\Delta_{k}&=&\sum_{0<\{m_{k}\frac{\theta}{\pi}\}<\delta}S^{-}_{M_{k}}(e^{\sqrt{-1}\theta})\label{3.32}
\eea
for every $k=1,\ldots,q$. Moreover we have
\bea
\min\left\{\ \left\{\frac{m_{k}\theta}{\pi}\right\},1-\left\{\frac{m_{k}\theta}{\pi}\right\}\right\}<\delta, \label{3.39a}
\eea
\bea
m_{k}\frac{\theta}{\pi}\in\Z, \;\ \rm if \;\ \frac{\theta}{\pi}\in\Q ,\label{3.40}
\eea
where $e^{\sqrt{-1}\theta}\in\sigma(M_{k})$, $\frac{\theta}{\pi}\in(0,2)$ and $\delta$ can be chosen as small as we want
(cf. (4.43)of\cite{LoZ1}).
\\\\ More precisely, by (4.10),(4.40), and(4.41) in \cite{LoZ1}, we
have
\bea
m_{k}=\left(\left[\frac{N}{M\hat{i}(\ga_{k},1)}\right]+\chi_{k}\right)M,  \ \ 1\leq k\leq q, \label{2.15}
\eea
where $\chi_{k}=0$ or $1$ for $1\leq k\leq q$
. Furthermore, given $M_{0}\in\N$, by the proof of Theorem 4.1 of \cite{LoZ1},
we may further require $M_{0}|N$ (since the closure of the set $\{\{N\upsilon\}: N\in \N, M_{0}|N\}$)
is still a closed additive subgroup of $T^{h}$ for some $h\in\N$, where we use the notations as
(4.21)-(4.22) in \cite{LoZ1}. Then we can use the step 2 in Theorem 4.1 of \cite{LoZ1} to get $N$).

In fact, by (4.40)-(4.41) of \cite{LoZ1}, let $\mu_{i}=\sum_{\theta\in(0,2\pi)}S^{-}_{M_{i}}(e^{\sqrt{-1}\theta})$
for $1\leq i\leq q$ and $\alpha_{i,j}=\frac{\theta_{j}}{\pi}$ where $e^{\sqrt{-1}\theta_{j}}\in \sigma (M_{i})$
for $1\leq j\leq \mu_{i}$ and $1\leq i\leq q$. As in (4.21) of \cite{LoZ1}, let $h=q+\sum_{1\leq i\leq q}\mu_{i}$ and
\bea
v=(\frac{1}{M\hat{i}(\gamma_{1},1)},\ldots,\frac{1}{M\hat{i}(\gamma_{q},1)},\frac{\alpha_{1,1}}{\hat{i}_{1}(\gamma_{1},1)},
\frac{\alpha_{1,2}}{\hat{i}_{1}(\gamma_{1},1)},\ldots \frac{\alpha_{1,\mu_{1}}}{\hat{i}_{1}(\gamma_{1},1)},
\frac{\alpha_{2,1}}{\hat{i}_{1}(\gamma_{2},1)},
\ldots,\frac{\alpha_{q,\mu_{q}}}{\hat{i}_{1}(\gamma_{q},1)}) .\label{3.35}
\eea
Then by (4.22) of \cite{LoZ1}, the above theorem is equivalent to find a vertex
\bea
\chi=(\chi_{1},\ldots,\chi_{q},\chi_{1,1},\chi_{1,2},\ldots,\chi_{1,\mu_{1}},\chi_{2,1},\ldots,\chi_{q,\mu_{q}})
\eea
of the cube $[0,1]^{h}$ and infinitely many integers $N\in\N$ such that
\bea
|\{Nv\}-\chi|<\varepsilon \label{3.36}
\eea
for any given $\varepsilon$ small enough.
\\\\
{\bf Theorem 2.6.}(cf. Theorem 4.2 of \cite{LoZ1}) Let $H$ be the closure of $\{\{mv\}|m\in \N\}$
in $T^{h}=(\R/Z)^{h}$ and $V=T_{0}\pi^{-1}H$ be the tangent space of $\pi^{-1}H$ at the origin in
$\R^{h}$, where $\pi: \R^{h}\rightarrow T^{h}$ is the projection map. Define
\bea
A(v)=V\setminus \cup_{v_{k}\in R\setminus Q}\{x=(x_{1},\ldots,x_{h})\in V|x_{k}=0\}.
\eea
Define $\psi(x)=0$ when $x\geq0$ and $\psi(x)=1$ when $x<0$. Then for any $a=(a_{1},\ldots,a_{h})
\in A(v)$, the vector
\bea
\chi(a)=(\psi(a_{1}),\ldots,\psi(a_{h})) \lb{3.36a}
\eea
makes (\ref{3.36}) hold for infinitely many $N\in \N$.
\\\\
Please note that when we choose $a\in V$  small enough, then $a+\chi(a)\in [0,1]^h$, this implies $(V+\chi(a))\cap[0,1]^h\neq\emptyset$,   and so we can require
 $N\in\N$ in (\ref{3.36}) satisfied   $\{Nv\}-\chi(a)\in V$.
 \\\\
{\bf Theorem 2.7.}(cf. Theorem 4.2 of \cite{LoZ1}) We have the following properties for $A(v)$:
\\
(i)When $v\in\R^{h}\setminus\Q^{h}$, then $dim V\geq1, 0\not\in A(v)\subset V$, $A(v)=-A(v)$
and $A(v)$ is open in $V$.
\\
(ii)When dim$V=1$, then $A(v)=V\setminus\{0\}$.
\\
(iii)When dim$V\geq 2$, $A(v)$ is obtained from $V$ by deleting all the coordinate hyperplanes
with dimension strictly smaller than dim$V$ from $V$.
\\\\
{\bf Remark  2.8.} In our choice of $(N,m_{1},\ldots,m_{q})$ in the proof of Theorem 2.5, we
can choose $M_{0}$ good enough such that $N\in\N$ further satisfies
\bea
\frac{N}{M\hat{i}(\gamma_{k},1)}\in \Z, \;\ {\rm for}\ \forall\  \hat{i}(\gamma_{k},1)\in \Q, \;\
k\in\{1,\ldots,q\}.\label{3.47}
\eea
Furthermore from (\ref{3.36}), we get $\hat{i}(\gamma_{k},1)\in \Q$ implies $\chi_{k}(a)=\psi(a_{k})=0$.
\\\\From the theorems above, we get a useful lemma below, this lemma is
very important in our proof of the main Theorem 2.\\\\
{\bf Lemma 2.9.} Let $v=(v_{1},v_{2},\ldots,v_{h})$ given by (\ref{3.35}). If $v_{i},v_{j}\in\R\setminus\Q
$ and $\frac{v_{j}}{v_{i}}=\frac{p}{q}\in Q^{+}(i<j)$, then for $\forall x\in V$, we have
$\frac{x_{j}}{x_{i}}=\frac{p}{q}$ and for $\forall a\in A(v)$, we have $\chi_{i}(a)=\chi_{j}(a)$. \\\\
{\bf Proof.} Since $H=\overline{\{\{mv\}|m\in\N\}}$, if $v_{i},v_{j}\in\R\setminus\Q
$ and $\frac{v_{j}}{v_{i}}=\frac{p}{q}\in Q^{+}(i<j)$, that means $v_{i}$, $v_{j}$
are rational dependent, then we restrict $H$ to the coordinate hyperplane
$$D=\{(0,\ldots,0,x_{i},0,\ldots,0,x_{j},0,\ldots,0)|x_{i},x_{j}\in R\}\subset R^{h},$$
we have the dynamical picture of $H\cap D$ below:
\begin{figure}[H]
\begin{minipage}[t]{0.333\linewidth}
\centering
\includegraphics[height=0.9\textwidth,width=1\textwidth]{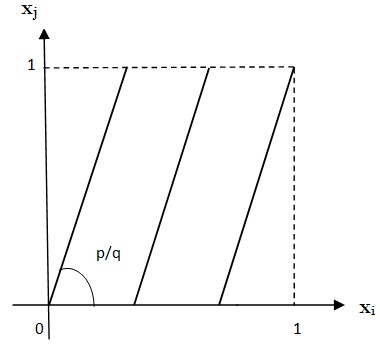}
  \caption{If $\frac{v_{j}}{v_{i}}=\frac{p}{q}>1$}
\label{fig:side:a}
\end{minipage}%
\begin{minipage}[t]{0.333\linewidth}
\centering
\includegraphics[height=0.9\textwidth,width=1\textwidth]{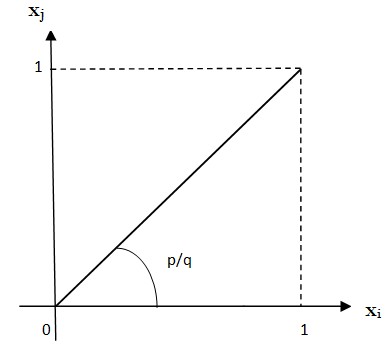}
     \caption{If $\frac{v_{j}}{v_{i}}=\frac{p}{q}=1$}
\label{fig:side:b}
\end{minipage}
\begin{minipage}[t]{0.333\linewidth}
\centering
\includegraphics[height=0.9\textwidth,width=1\textwidth]{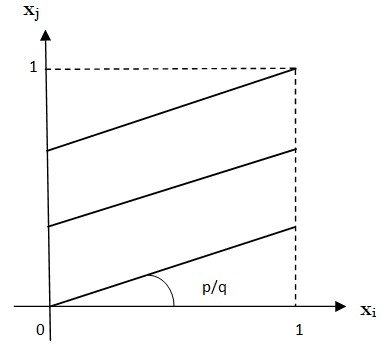}
     \caption{If $\frac{v_{j}}{v_{i}}=\frac{p}{q}<1$}
\label{fig:side:c}
\end{minipage}
\end{figure}
\noindent hence for any $x\in V=T_{0}\pi^{-1}H, x_{i}\neq0$, we have $\frac{x_{j}}{x_{i}}=\frac{p}{q}\in\Q^{+}$. In
particular, if $a\in A(v)$, then $a_{i}>0, a_{j}>0$ or $a_{i}<0, a_{j}<0$, hence $\chi_{i}(a)=\chi_{j}(a)=0$
or $\chi_{i}(a)=\chi_{j}(a)=1$, this completes the proof.  $\Box$

Furthermore, as Definition 1.1 of \cite{LoZ1}, we define
\\\\
{\bf Definition 2.10.} For $\alpha\in(1,2)$, we define a map $\varrho_{n}(\Sigma):\H(2n)\rightarrow \N\cup\{+\infty\}$
\bea
\varrho_{n}(\Sigma)=\left\{\matrix{+\infty, &\quad {\rm if}\;\; ^{\#}\mathcal{V}(\Sigma,\alpha)=+\infty, \cr
\min\left\{\left[\frac{i(x,1)+2S^{+}(x)-\nu(x,1)+n}{2}\right]
\mid (\tau,x)\in \mathcal{V}(\Sigma,\alpha)\right\}, &\quad {\rm if}\;\; ^{\#}\mathcal{V}(\Sigma,\alpha)<+\infty. \cr}\right.
\lb{6.17}
\eea
\\\\
As in \cite{LoZ1},  we denote  the elements of $\mathcal{V}_{\infty}(\Sigma,\alpha)$ by
$$ \mathcal{V}_{\infty}(\Sigma,\alpha)=\{ [(\tau_j,x_j)] | j=1,\cdots,q   \},  $$
where  $(\tau_j,x_j)\in \J(\Sigma,\alpha)$ for $j=1,\cdots,q$.
\\\\
{\bf Theorem 2.11.} For given $a\in A(v)$, we define $\chi\equiv\chi(a)=(\psi(a_{1}),\ldots,\psi(a_{h}))$
by (\ref{3.36a}). Let $(N,m_{1},\ldots,m_{q})\in\N^{q+1}$ be given in Remark 2.8. Then by the proof Theorem 5.1
in \cite{LoZ1}, for each $s=1,\ldots,\varrho_{n}(\Sigma)$, there exists a unique $j(s)\in\{1,\ldots,q\}$
and an injection map $p: N+K\rightarrow \mathcal{V}_{\infty}(\Sigma,\alpha)\times \N$ such that
$p(N-s+1)=([(\tau_{j(s)},x_{j(s)})],2m_{j(s)})$ and
\bea
i(x_{j(s)},2m_{j(s)})\leq2N-2s+n\leq i(x_{j(s)},2m_{j(s)})+\nu(x_{j(s)},2m_{j(s)})-1, \label{3.42}
\eea
Then for any $s_{1}, s_{2}\in\{1,\ldots,\varrho_{n}(\Sigma)\}$
with $s_{1}<s_{2}$, we have:
\bea
\left(\left[\frac{N}{MD_{j(s_{2})}}\right]+\chi_{j(s_{2})}(a)\right)MD_{j(s_{2})}
&<&\left(\left[\frac{N}{MD_{j(s_{1})}}\right]+\chi_{j(s_{1})}(a)\right)MD_{j(s_{1})} \label{3.43}
\eea
and
\bea
i(x_{j(s)},2m_{j(s)})&=&2(N+\Delta_{j(s)})-(S_{M_{j(s)}}^{+}(1)+C(M_{j(s)})),\label{3.44} \\
2s&\geq& n+S_{M_{j(s)}}^{+}(1)+C(M_{j(s)})-2\Delta_{j(s)}-\nu(x_{j(s)},2m_{j(s)})+1,\label{3.45}\\
2s&\leq& n+S_{M_{j(s)}}^{+}(1)+C(M_{j(s)})-2\Delta_{j(s)}, \label{3.52}
\eea
where $D_{j(s)}=\hat{i}(x_{j(s)},1)$.
\\\\
{\bf Proof.}
Since $s_{1}<s_{2}$, we have $N-s_{1}+1>N-s_{2}+1$, then from Theorem 2.3, we have
\bea \nonumber
\hat{i}(x_{j(s_{2})},2m_{j(s_{2})})<\hat{i}(x_{j(s_{1})},2m_{j(s_{1})}),
\eea
the property of the mean index $\hat{i}(x,m)=m \hat{i}(x,1)$ implies that
\bea
2m_{j(s_{2})}\hat{i}(x_{j(s_{2})},1)<2m_{j(s_{1})}\hat{i}(x_{j(s_{1})},1)).
\eea
From the definition of $m_{j(s)}=([\frac{N}{M\hat{i}(\gamma_{j(s)},1)}]+\chi_{j(s)}(a))M$ and
$D_{j(s)}=\hat{i}(x_{j(s)},1)$
we get (\ref{3.43}). In order to prove formula (\ref{3.44}), we need some identities
(\ref{3.30}), (\ref{3.31}) and (\ref{3.23}) below
\bea \nonumber
I(j(s),m_{j(s)})&=&N+\Delta_{j(s)},
\eea where
\bea \nonumber
I(j(s),m_{j(s)})&=&m_{j(s)}(i(\gamma_{j(s)},1)+S_{M_{j(s)}}^{+}(1)-C(M_{j(s)}))\\ \nonumber
&+&\sum_{\theta\in(0,2\pi)}E(\frac{m_{j(s)}\theta}{\pi})
S_{M_{j(s)}}^{-}(e^{\sqrt{-1}\theta})),
\eea
\bea\nonumber
i(\gamma_{j(s)},m_{j(s)})&=&m_{j(s)}(i(\gamma_{j(s)},1)+S^{+}_{M_{j(s)}}(1)-C(M_{j(s)})) \\ \nonumber
&+&2\sum_{\theta\in(0,2\pi)}E(\frac{m_{j(s)}\theta}{2\pi})S^{-}_{M_{j(s)}}(e^{\sqrt{-1}\theta})
-(S_{M_{j(s)}}^{+}(1)+C(M_{j(s)})).
\eea
where $C(M_{j(s)})=\sum_{0<\theta<2\pi}S^{-}_{M_{j(s)}}(e^{\sqrt{-1}\theta})$.

Simple calculations show that
\bea \nonumber
i(\gamma_{j(s)},2m_{j(s)})&=&2I(j(s),m_{j(s)})-(S_{M_{j(s)}}^{+}(1)+C(M_{j(s)}) \\ \nonumber
&=&2(N+\Delta_{j(s)})-(S_{M_{j(s)}}^{+}(1)+C(M_{j(s)}).
\eea
hence we get formula (\ref{3.44}). On the other hand, it's easy to show that formulas (\ref{3.42}),(\ref{3.44})
imply (\ref{3.45}), (\ref{3.52}). $\Box$
\\\\
{\bf Corollary 2.12.} Further properties of inequalities (\ref{3.43}):

i) If $D_{j(s_{2})}\in\Q$, then $\chi_{j(s_{2})}(a)=0$ and
\bea \nonumber
N=\left(\left[\frac{N}{MD_{j(s_{2})}}\right]+\chi_{j(s_{2})}(a)\right)MD_{j(s_{2})}
&<&\left(\left[\frac{N}{MD_{j(s_{1})}}\right]+\chi_{j(s_{1})}(a)\right)MD_{j(s_{1})}
\eea
with $D_{j(s_{1})}\in\R\setminus\Q$, $\chi_{j(s_{1})}(a)=1$.

ii) If $\chi_{j(s_{2})}(a)=1$, then $D_{j(s_{2})}\in\R\setminus\Q$ and
\bea \nonumber
N<\left(\left[\frac{N}{MD_{j(s_{2})}}\right]+\chi_{j(s_{2})}(a)\right)MD_{j(s_{2})}
&<&\left(\left[\frac{N}{MD_{j(s_{1})}}\right]+\chi_{j(s_{1})}(a)\right)MD_{j(s_{1})}
\eea
with $D_{j(s_{1})}\in\R\setminus\Q$, $\chi_{j(s_{1})}(a)=1$.

iii)  \bea  \chi_{j(s_{2})}(a)\leq \chi_{j(s_{1})}(a).\label{ad1} \eea
\\\\
{\bf Proof.} From Remark 2.8, we know that $D_{j(s_{2})}\in\Q$ implies $\frac{N}{MD_{j(s_{2})}}\in\Z$, then
(\ref{3.36}) implies that $\chi_{j(s_{2})}(a)=0$, hence the formula
$([\frac{N}{MD_{j(s_{2})}}]+\chi_{j(s_{2})}(a))MD_{j(s_{2})}=N$. For this case,
it's easy to check that inequality (\ref{3.43}) holds if and only if $D_{j(s_{1})}\in\R\setminus\Q$, $\chi_{j(s_{1})}(a)=1$.
This completes the proof of i). From i) we know if $D_{j(s_{2})}\in\Q$
then $\chi_{j(s_{2})}(a)=0$, so $\chi_{j(s_{2})}(a)=1$ implies that $D_{j(s_{2})}\in\R\setminus\Q$, easy computation shows that
\bea
N=\left(\frac{N}{MD_{j(s_{2})}}\right)MD_{j(s_{2})}
&<&\left(\left[\frac{N}{MD_{j(s_{2})}}\right]+\chi_{j(s_{2})}(a)\right)MD_{j(s_{2})}.
\eea
This combine with inequality (\ref{3.43}), we get ii).  To prove iii), please note that if  $\chi_{j(s_{2})}(a)=0$, then  (\ref{ad1}) is obviously right,  the case  $\chi_{j(s_{2})}(a)=1$ is from ii).
\\\\
{\bf Remark 2.13}  It is proved in \cite{LoZ1} that for $s=1,...,\varrho_{n}(\Sigma)$,  $p(N-s+1)$ are geometric different, thus there are at least $\varrho_{n}(\Sigma)$ closed characteristics, and moreover at least $\varrho_{n}(\Sigma)-1$ among them have irrational mean index and
$p(N)$ is elliptic.

\setcounter{equation}{0}
\section{ Proofs of the Theorems 1.1 and 1.2}

In this section, we prove Theorems 1.1 and 1.2 based on the
index iteration theory developed by Y. Long and his coworkers.
Some notations for the Maslov-type index can be found in Section 4.  The basic normal form $R(\theta_{j})$
($N_{2}(\omega_{j},u_{j})
$; $N_{2}(\lambda_{j},\nu_{j}))$ given in Theorem 4.7 is called rational normal form
, if $\frac{\theta_{j}}{\pi}\in\Q (\frac{\alpha_{j}}{\pi}\in\Q; \frac{\beta_{j}}{\pi}\in\Q)$.
\\\\
{\bf Lemma 3.1.} For any fix $a\in A(v)$, let inject map
$p(N-s+1)=([(\tau_{j(s)},x_{j(s)})],2m_{j(s)}), s\in\{1,\ldots,\varrho_{n}(\Sigma)\}$
given in Theorem 2.11, if $x_{j(2)}$ is not an elliptic closed characteristic, then $\chi_{j(2)}(a)=0$
implies that $\hat{i}(x_{j(2)},1)\in\Q$.
\\\\
{\bf Proof.} From Theorem 4.7, we have the symplectic decomposition
\bea
\gamma_{j(2)}&\simeq&N_{1}(1,1)^{\diamond p_{-}}\diamond I_{2p_{0}} \diamond N_{1}(1,-1)^{\diamond p_{+}} \nonumber
\diamond N_{1}(-1,1)^{\diamond q_{-}}\diamond -I_{2q_{0}} \diamond N_{1}(-1,-1)^{\diamond q_{+}} \\ \nonumber
&\diamond& R(\theta_{1})\diamond\cdots R(\theta_{r}) \diamond N_{2}(\omega_{1},u_{1})\diamond
\cdots \diamond N_{2}(\omega_{r_{*}},u_{r_{*}}) \\
&\diamond& N_{2}(\lambda_{1},\upsilon_{1}) \diamond
\cdots \diamond N_{2}(\lambda_{r_{0}},\upsilon_{r_{0}})\diamond M_{k},
\eea
for this decomposition, the number of the rational normal form
in $\{R(\theta_{1}),\ldots,R(\theta_{r})\}$ is denoted by $\tilde{r}$. Similarly,
for set $\{N_{2}(\omega_{1},u_{1}),\ldots,N_{2}(\omega_{r_{*}},u_{r_{*}})\}$
and $\{N_{2}(\lambda_{1},\nu_{1}),\ldots,N_{2}(\lambda_{r_{0}},\nu_{r_{0}})\}$,
the number of rational normal form is
denoted by $\tilde{r}_{*}$ and $\tilde{r}_{0}$ respectively, then from (\ref{3.27}) we have a further
estimation of the variable $\nu(\gamma_{j(2)},2m_{j(2)})$ below
\bea \nonumber
\nu(\gamma_{j(2)},2m_{j(2)})&=&\nu(\gamma_{j(2)},1)+q_{-}+2q_{0}+q_{+}+2(r+r_{*}+r_{0})\nonumber\\
&&-2(r-\tilde{r}+r_{*}-\tilde{r}_{*}+r_{0}-\tilde{r}_{0})\nonumber \\
&=&p_{-}+2p_{0}+p_{+}+q_{-}+2q_{0}+q_{+}+2(\tilde{r}+\tilde{r}_{*}+\tilde{r}_{0}).\label{3.2}
\eea
Now we proof the lemma by contradiction. Assume that $\chi_{j(2)}(a)=0$
and $\hat{i}(x_{j(2)},1)\in\R\setminus\Q$, then (\ref{3.28}) implies that
at least one of $\frac{\theta_{1}}{\pi}, \frac{\theta_{2}}{\pi},\ldots,\frac{\theta_{r}}{\pi}$
is irrational number, hence $r-\tilde{r}\geq1$ and
\bea \nonumber
\left\{m_{j(2)}D_{j(2)}\right\}&=&
\left\{m_{j(2)}\left(i(x_{j(2)},1)+p_{-}+p_{0}-r+\sum_{j=1}^{r}\frac{\theta_{j}}{\pi}\right)\right\}\\ \nonumber
&=&\left\{m_{j(2)}\sum_{\frac{\theta_{j}}{\pi}\in \R\setminus\Q}\frac{\theta_{j}}{\pi}\right\}\\
&\leq&\sum_{\frac{\theta_{j}}{\pi}\in \R\setminus\Q}\left\{m_{j(2)}\frac{\theta_{j}}{\pi}\right\}, \label{4.7}
\eea
where $m_{j(2)}=([\frac{N}{MD_{j(2)}}]+\chi_{j(2)}(a))M=[\frac{N}{MD_{j(2)}}]M$. The second equality comes from
(\ref{3.40}) in Theorem 2.5.

On the other hand,
\bea \nonumber
\left\{m_{j(2)}D_{j(2)}\right\}&=&
\left\{\left[\frac{N}{MD_{j(2)}}\right]MD_{j(2)}\right\}\\
&=&\left\{N-\left\{\frac{N}{MD_{j(2)}}\right\}MD_{j(2)}\right\}, \label{4.8}
\eea
and from (\ref{3.35}) and (\ref{3.36}), we get that $\{\frac{N}{MD_{j(2)}}\}=
\left|\{\frac{N}{MD_{j(2)}}\}-\chi_{j(2)}(a)\right|<\varepsilon$ ($\chi_{j(2)}(a)=0$) for any given
$\varepsilon$ small enough, let $\varepsilon<\frac{1-\delta}{MD_{j(2)}}$,
, where $\delta$ is given in Theorem 2.5,
hence $\{m_{j(2)}D_{j(2)}\}>\delta,$
this combines with (\ref{4.7}), we have
\bea
\delta
<\sum_{\frac{\theta_{j}}{\pi}\in \R\setminus\Q}\left\{m_{j(2)}\frac{\theta_{j}}{\pi}\right\},
\eea
hence at least one of the elements in $\{\frac{\theta_{j}}{\pi}|\frac{\theta_{j}}{\pi}\in \R\setminus\Q, j=1,\ldots,r\}$
satisfies $\{m_{j(2)}\frac{\theta_{j}}{\pi}\}\not\in(0,\delta)$. We have the estimation of the variable $\Delta_{j(2)}$ below
\bea \nonumber
\Delta_{j(2)}&=&\sum_{0<\{m_{j(2)}\frac{\theta}{\pi}\}<\delta}S_{M_{j(2)}}^{-}(e^{\sqrt{-1}\theta})\nonumber\\
&\leq&\sum_{\frac{\theta}{\pi}\in\R\setminus\Q}S_{M_{j(2)}}^{-}(e^{\sqrt{-1}\theta})-1\nonumber\\
&=& r-\tilde{r}-1+2(r_{*}-\tilde{r}_{*}),\label{4.11}
\eea
where  the last equality comes from the calculation of the splitting number of the basic normal given in Definition 4.2.
In order to prove the Lemma, we rewrite the useful inequality (\ref{3.45}) and equalities
(\ref{3.25}), (\ref{3.29}), (\ref{4.34}), (\ref{3.2}) below
\bea \nonumber
2s&\geq& n+S_{M_{j(s)}}^{+}(1)+C(M_{j(s)})-2\Delta_{j(s)}-\nu(x_{j(s)},2m_{j(s)})+1 ,\;\ (s=2)\nonumber\\
n&=&p_{-}+p_{0}+p_{+}+q_{-}+q_{0}+q_{+}+r+2r_{*}+2r_{0}+k\nonumber\\
S_{M_{j(2)}}^{+}(1)&=&p_{-}+p_{0}\nonumber\\
C(M_{j(2)})&=&\sum_{0<\theta<2\pi}S_{M_{j(2)}}^{-}(e^{\sqrt{-1}\theta})=q_{0}+q_{+}+r+2r_{*}\nonumber\\
\nu(\gamma_{j(2)},2m_{j(2)})
&=&p_{-}+2p_{0}+p_{+}+q_{-}+2q_{0}+q_{+}+2(\tilde{r}+\tilde{r}_{*}+\tilde{r}_{0})\nonumber
\eea
This combine with inequality (\ref{4.11}), easy computation shows that for $s=2$, we have
\bea
4\geq p_{-}+q_{+}+2(r_{0}-\tilde{r}_{0})+2\tilde{r}_{*}+k+3
\eea
From Remark 2.4, we always have $p_{-}\geq1$. On the other hand, from the condition of the lemma,
we know $x_{j(2)}$ is not an elliptic closed characteristic, this implies that $k\geq1$, hence we
get $4\geq5.$ This contradiction completes the proof.  $\Box$
\\\\
{\bf Corollary 3.2.} For any fix $a\in A(v)$, let inject map
$p(N-s+1)=([(\tau_{j(s)},x_{j(s)})],2m_{j(s)}), s\in\{1,\ldots,\varrho_{n}(\Sigma)\}$
given in Theorem 2.11. If $x_{j(2)}$ is not elliptic, then
$\chi_{j(1)}(a)=1$ and $\hat{i}(x_{j(1)},1)\in\R\setminus\Q$.\\\\
{\bf Proof.} If $\chi_{j(2)}(a)=1$, from ii) of Corollary 2.12 we have
$\chi_{j(1)}(a)=1$ and $\hat{i}(x_{j(1)},1)\in\R\setminus\Q$. If $\chi_{j(2)}(a)=0$,
then from Lemma 3.1, we have $\hat{i}(x_{j(2)},1)\in\Q$. This combines with i) of Corollary 2.12 we get
$\hat{i}(x_{j(1)},1)\in\R\setminus\Q$ and $\chi_{j(1)}(a)=1$.  $\Box$

Now we start to proof Theorem 1.1 and Theorem 1.2.
\\\\
{\bf Proof of Theorem 1.1.} For any fix $a\in A(v)$, let inject map
$p(N-s+1)=([(\tau_{j(s)},x_{j(s)})],\\ 2m_{j(s)}),s\in\{1,\ldots,\varrho_{n}(\Sigma)\},$
given in Theorem 2.11. From Remark 2.13, we know $x_{j(1)}$ is elliptic.
If $x_{j(2)}$ is also elliptic, then the proof is complete.
If $x_{j(2)}$ is not an elliptic closed characteristic, from Corollary 3.2, we get
$\chi_{j(1)}(a)=1$ and $\hat{i}(x_{j(1)},1)\in\R\setminus\Q$. Now from Theorem 2.7, we
can choose $-a\in A(v)$, then Theorem 2.11 says that for $-a\in A(v)$,
we still have $(\tilde{N},\tilde{m}_{1},\ldots,\tilde{m}_{q})$
, $\tilde{j}(s)$ and inject map
$p(\tilde{N}-s+1)=([(\tau_{\tilde{j}(s)},x_{\tilde{j}(s)})],2\tilde{m}_{\tilde{j}(s)}), s\in\{1,\ldots,\varrho_{n}(\Sigma)\}$.
If $\tilde{j}(1)\neq j(1)$, then $x_{j(1)}, x_{\tilde{j}(1)}$ are two different elliptic closed
characteristics, then we complete the proof. If $\tilde{j}(1)= j(1)$, from the definition of
$\chi(a)$, we know that $\chi_{j(1)}(a)=1$ implies $\chi_{j(1)}(-a)=0$, hence
$\chi_{\tilde{j}(1)}(-a)=\chi_{j(1)}(-a)=0$, this combine with i), iii) of Corollary 2.12,
we have $\chi_{\tilde{j}(2)}(-a)=0$ and $\hat{i}(x_{\tilde{j}(2)},1)\in\R\setminus\Q$, but
Lemma 3.1 still holds in the case
$-a\in A(v)$, that means if $x_{\tilde{j}(2)}$ is not elliptic, we should have
$\chi_{\tilde{j}(2)}(-a)=0$ implies $\hat{i}(x_{\tilde{j}(2)},1)\in\Q$, this contradiction completes
the proof.    $\Box$\\\\
{\bf Proof of Theorem 1.2.} For the $\varrho_{n}(\Sigma)$ geometrically
distinct closed characteristics in Remark 2.13,  for any two closed characteristics$[(\tau,x)], [(\tilde{\tau},\tilde{x})]$, we know that if one of
them
has rational mean index, then anther must has irrational mean index, hence
for this case, the theorem is true. Now we can assume that $\hat{i}(x,1), \hat{i}(\tilde{x},1)\in\R\setminus\Q$.
 For this case, we get the  proof by contradiction.

Assume $\frac{\hat{i}(\tilde{x},1)}{\hat{i}(x,1)}=\frac{p}{q}\in Q^{+}$, then from
Theorem 2.11, we know there exist $s_{1}, s_{2}\in \{1,2,\ldots,q\}$ such that
$x=x_{j(s_{1})}$ and $\tilde{x}=x_{j(s_{2})}$. Without loss of generality, we can assume that
$s_{1}<s_{2}$, then we have (\ref{3.43})
\bea \nonumber
\left(\left[\frac{N}{MD_{j(s_{2})}}\right]+\chi_{j(s_{2})}(a)\right)MD_{j(s_{2})}
&<&\left(\left[\frac{N}{MD_{j(s_{1})}}\right]+\chi_{j(s_{1})}(a)\right)MD_{j(s_{1})},
\eea
hence we have
\bea
\left(\left\{\frac{N}{MD_{j(s_{2})}}\right\}-\chi_{j(s_{2})}(a)\right)MD_{j(s_{2})}
&>&\left(\left\{\frac{N}{MD_{j(s_{1})}}\right\}-\chi_{j(s_{1})}(a)\right)MD_{j(s_{1})}, \label{4.2}
\eea
where $D_{j(s)}=\hat{i}(x_{j(s)},1)$. Let
\bea \nonumber
v=(\frac{1}{M\hat{i}(\gamma_{1},1)},\ldots,\frac{1}{M\hat{i}(\gamma_{q},1)},\frac{\alpha_{1,1}}{\hat{i}_{1}(\gamma_{1},1)},
\frac{\alpha_{1,2}}{\hat{i}_{1}(\gamma_{1},1)},\ldots \frac{\alpha_{1,\mu_{1}}}{\hat{i}_{1}(\gamma_{1},1)},
\frac{\alpha_{2,1}}{\hat{i}_{1}(\gamma_{2},1)},
\ldots,\frac{\alpha_{q,\mu_{q}}}{\hat{i}_{1}(\gamma_{q},1)})
\eea
given by (\ref{3.35}), where $\gamma_{k}$ is the associated symplectic path of
$[(\tau_{k},x_{k})]\in\mathcal{V}_{\infty}(\Sigma, \alpha)$, then we have
\bea
\frac{v_{j(s_{2})}}{v_{j(s_{1})}}=\frac{\frac{1}{M\hat{i}(\gamma_{j(s_{2})},1)}}{\frac{1}{M\hat{i}(\gamma_{j(s_{1})}£¬1)}}
=\frac{\hat{i}(x,1)}{\hat{i}(\tilde{x},1)}=\frac{q}{p}. \label{4.3}
\eea
On the other hand, from Theorem 2.6,   for fixed $a\in V$,   we choose  $N\in\N$ such that  $\{Nv\}-\chi(a)$ small enough,  recall that we also have $\{Nv\}-\chi(a)\in V$.
From  Lemma 2.9,  we   we get
\bea
\frac{\{Nv_{j(s_{2})}\}-\chi_{j(s_{2})}(a)}{\{Nv_{j(s_{1})}\}-\chi_{j(s_{1})}(a)}
=\frac{\{\frac{N}{MD_{j(s_{2})}}\}-\chi_{j(s_{2})}(a)}{\{\frac{N}{MD_{j(s_{1})}}\}-\chi_{j(s_{1})}(a)}
=\frac{q}{p}. \label{4.4}
\eea
(\ref{4.3}), (\ref{4.4}) imply
\bea
\left(\left\{\frac{N}{MD_{j(s_{1})}}\right\}-\chi_{j(s_{1})}(a)\right)MD_{j(s_{1})}
&=&\left(\left\{\frac{N}{MD_{j(s_{2})}}\right\}-\chi_{j(s_{2})}(a)\right)MD_{j(s_{2})}
\eea
This contradiction with (\ref{4.2}), then the proof is complete.   $\Box$\\

\setcounter{equation}{0}
\section{Appendix:  index iteration theory for closed characteristics}

In this section, we  briefly  recall the  index theory for symplectic paths.  The index theory  is introduced  by Conley and Zehnder  \cite{CoZ1} and developed by Long and others (see \cite{Lon5} for details).

As usual, the symplectic group $\Sp(2n)$ is defined by
$$ \Sp(2n) = \{M\in {\rm GL}(2n,\R)\,|\,M^TJM=J\}, $$
whose topology is induced from that of $\R^{4n^2}$. For $\tau>0$ we are interested
in paths in $\Sp(2n)$:
$$ \P_{\tau}(2n) = \{\ga\in C([0,\tau],\Sp(2n))\,|\,\ga(0)=I_{2n}\}, $$
which is equipped with the topology induced from that of $\Sp(2n)$. The
following real function was introduced in \cite{Lon4}:
$$ D_{\om}(M) = (-1)^{n-1}\ol{\om}^n\det(M-\om I_{2n}), \qquad
          \forall \om\in\U,\, M\in\Sp(2n). $$
Thus for any $\om\in\U$ the following codimension $1$ hypersurface in $\Sp(2n)$ is
defined in \cite{Lon4}:
$$ \Sp(2n)_{\om}^0 = \{M\in\Sp(2n)\,|\, D_{\om}(M)=0\}.  $$
For any $M\in \Sp(2n)_{\om}^0$, we define a co-orientation of $\Sp(2n)_{\om}^0$
at $M$ by the positive direction $\frac{d}{dt}Me^{t\ep J}|_{t=0}$ of
the path $Me^{t\ep J}$ with $0\le t\le 1$ and $\ep>0$ being sufficiently
small. Let
\bea
\Sp(2n)_{\om}^{\ast} &=& \Sp(2n)\bs \Sp(2n)_{\om}^0,   \nn\\
\P_{\tau,\om}^{\ast}(2n) &=&
      \{\ga\in\P_{\tau}(2n)\,|\,\ga(\tau)\in\Sp(2n)_{\om}^{\ast}\}, \nn\\
\P_{\tau,\om}^0(2n) &=& \P_{\tau}(2n)\bs  \P_{\tau,\om}^{\ast}(2n).  \nn\eea
For any two continuous arcs $\xi$ and $\eta:[0,\tau]\to\Sp(2n)$ with
$\xi(\tau)=\eta(0)$, it is defined as usual:
$$ \eta\ast\xi(t) = \left\{\matrix{
            \xi(2t), & \quad {\rm if}\;0\le t\le \tau/2, \cr
            \eta(2t-\tau), & \quad {\rm if}\; \tau/2\le t\le \tau. \cr}\right. $$
Given any two $2m_k\times 2m_k$ matrices of square block form
$M_k=\left(\matrix{A_k&B_k\cr
                                C_k&D_k\cr}\right)$ with $k=1, 2$,
as in \cite{Lon5}, the $\;\dm$-product of $M_1$ and $M_2$ is defined by
the following $2(m_1+m_2)\times 2(m_1+m_2)$ matrix $M_1\dm M_2$:
$$ M_1\dm M_2=\left(\matrix{A_1&  0&B_1&  0\cr
                               0&A_2&  0&B_2\cr
                             C_1&  0&D_1&  0\cr
                               0&C_2&  0&D_2\cr}\right). \nn$$  
Denote by $M^{\dm k}$ the $k$-fold $\dm$-product $M\dm\cdots\dm M$. Note
that the $\dm$-product of any two symplectic matrices is symplectic. For any two
paths $\ga_j\in\P_{\tau}(2n_j)$ with $j=0$ and $1$, let
$\ga_0\dm\ga_1(t)= \ga_0(t)\dm\ga_1(t)$ for all $t\in [0,\tau]$.

A special path $\xi_n\in\P_{\tau}(2n)$ is defined by
\be \xi_n(t) = \left(\matrix{2-\frac{t}{\tau} & 0 \cr
                                             0 &  (2-\frac{t}{\tau})^{-1}\cr}\right)^{\dm n}
         \qquad {\rm for}\;0\le t\le \tau.  \lb{6.1}\ee
  \\\\
{\bf Definition 4.1.} (cf. \cite{Lon4}, \cite{Lon5}) {\it For any $\om\in\U$ and
$M\in \Sp(2n)$, define
\be  \nu_{\om}(M)=\dim_{\C}\ker_{\C}(M - \om I_{2n}).  \lb{6.2}\ee
For any $\tau>0$ and $\ga\in \P_{\tau}(2n)$, define
\be  \nu_{\om}(\ga)= \nu_{\om}(\ga(\tau)).  \lb{6.3}\ee

If $\ga\in\P_{\tau,\om}^{\ast}(2n)$, define
\be i_{\om}(\ga) = [\Sp(2n)_{\om}^0: \ga\ast\xi_n],  \lb{6.4}\ee
where the right hand side of (\ref{6.4}) is the usual homotopy intersection
number, and the orientation of $\ga\ast\xi_n$ is its positive time direction under
homotopy with fixed end points.

If $\ga\in\P_{\tau,\om}^0(2n)$, we let $\mathcal{F}(\ga)$
be the set of all open neighborhoods of $\ga$ in $\P_{\tau}(2n)$, and define
\be i_{\om}(\ga) = \sup_{U\in\mathcal{F}(\ga)}\inf\{i_{\om}(\beta)\,|\,
                       \beta\in U\cap\P_{\tau,\om}^{\ast}(2n)\}.
               \lb{6.5}\ee
Then
$$ (i_{\om}(\ga), \nu_{\om}(\ga)) \in \Z\times \{0,1,\ldots,2n\}, $$
is called the index function of $\ga$ at $\om$. }

For any symplectic path $\ga\in\P_{\tau}(2n)$ and $m\in\N$,  we
define its $m$-th iteration $\ga^m:[0,m\tau]\to\Sp(2n)$ by
\be \ga^m(t) = \ga(t-j\tau)\ga(\tau)^j, \qquad
  {\rm for}\quad j\tau\leq t\leq (j+1)\tau,\;j=0,1,\ldots,m-1.
     \lb{6.6}\ee
We still denote the extended path on $[0,+\infty)$ by $\ga$.
\\\\
{\bf Definition 4.2} (cf. \cite{Lon4}, \cite{Lon5}) {\it For any $\ga\in\P_{\tau}(2n)$,
we define
\be (i(\ga,m), \nu(\ga,m)) = (i_1(\ga^m), \nu_1(\ga^m)), \qquad \forall m\in\N.
   \lb{6.7}\ee
The mean index $\hat{i}(\ga,m)$ per $m\tau$ for $m\in\N$ is defined by
\be \hat{i}(\ga,m) = \lim_{k\to +\infty}\frac{i(\ga,mk)}{k}. \lb{6.8}\ee
For any $M\in\Sp(2n)$ and $\om\in\U$, the {\it splitting numbers} $S_M^{\pm}(\om)$
of $M$ at $\om$ are defined by
\be S_M^{\pm}(\om)
     = \lim_{\ep\to 0^+}i_{\om\exp(\pm\sqrt{-1}\ep)}(\ga) - i_{\om}(\ga),
   \lb{6.9}\ee
for any path $\ga\in\P_{\tau}(2n)$ satisfying $\ga(\tau)=M$.}

For $\Sigma\in\H(2n)$ and $\alpha\in(1,2)$, let $(\tau,x)\in\J(\Sigma,\alpha)$. we define
\bea
S^{+}(x)&=&S^{+}_{\gamma_{x}(\tau)}(1),\\
(i(x,m),\nu(x,m))&=&(i(\gamma_{x},m),\nu(\gamma_{x},m)),\\
\hat{i}(x,m)&=&\hat{i}(\gamma_{x},m),
\eea
For all $m\in\N$, where $\gamma_{x}$ is the associated symplectic path of $(\tau,x)$.

For a given path $\gamma\in {\cal P}_{\tau}(2n)$ we consider to deform
it to a new path $\eta$ in ${\cal P}_{\tau}(2n)$ so that
\begin{equation}
i_1(\gamma^m)=i_1(\eta^m),\quad \nu_1(\gamma^m)=\nu_1(\eta^m), \quad
         \forall m\in {\bf N}, \label{6.10}
\end{equation}
and that $(i_1(\eta^m),\nu_1(\eta^m))$ is easy enough to compute. This
leads to finding homotopies $\delta:[0,1]\times[0,\tau]\to {\rm Sp}(2n)$
starting from $\gamma$ in ${\cal P}_{\tau}(2n)$ and keeping the end
points of the homotopy always stay in a certain suitably chosen maximal
subset of ${\rm Sp}(2n)$ so that (\ref{6.10}) always holds. In fact,  this
set was first discovered in \cite{Lon2} as the path connected component
$\Omega^0(M)$ containing $M=\gamma(\tau)$ of the set
\begin{eqnarray}
  \Omega(M)=\{N\in{\rm Sp}(2n)\,&|&\,\sigma(N)\cap{\bf U}=\sigma(M)\cap{\bf U}\;
{\rm and}\;  \nonumber\\
 &&\qquad \nu_{\lambda}(N)=\nu_{\lambda}(M), \;\forall\,
\lambda\in\sigma(M)\cap{\bf U}\}. \label{6.11}
\end{eqnarray}
Here $\Omega^0(M)$ is called the {\it homotopy component} of $M$ in
${\rm Sp}(2n)$.

In \cite{Lon2}-\cite{Lon5}, the following symplectic matrices were introduced
as {\it basic normal forms}:
\begin{eqnarray}
D(\lambda)=\left(\matrix{\lm & 0\cr
         0  & \lm^{-1}\cr}\right), &\quad& \lm=\pm 2,\lb{6.12}\\
N_1(\lm,b) = \left(\matrix{\lm & b\cr
         0  & \lm\cr}\right), &\quad& \lm=\pm 1, b=\pm1, 0, \lb{6.13}\\
R(\th)=\left(\matrix{\cos\th & -\sin\th\cr
        \sin\th  & \cos\th\cr}\right), &\quad& \th\in (0,\pi)\cup(\pi,2\pi),
                     \lb{6.14}\\
N_2(\om,b)= \left(\matrix{R(\th) & b\cr
              0 & R(\th)\cr}\right), &\quad& \th\in (0,\pi)\cup(\pi,2\pi),
                     \lb{6.15}\end{eqnarray}
where $b=\left(\matrix{b_1 & b_2\cr
               b_3 & b_4\cr}\right)$ with  $b_i\in\R$ and  $b_2\not=b_3$, $\omega=e^{\sqrt{-1}\theta}$.
We call $N_{2}(\om,b)$ is nontrivial if $(b_{2}-b_{3})\sin \theta<0$ and
$N_{2}(\om,b)$ is trivial if $(b_{2}-b_{3})\sin \theta>0$.

Splitting numbers possess the following properties:
\\\\
{\bf Lemma 4.3.} (cf. \cite{Lon2} and Lemma 9.1.5 of \cite{Lon5}) {\it Splitting
numbers $S_M^{\pm}(\om)$ are well defined, i.e., they are independent of the choice
of the path $\ga\in\P_\tau(2n)$ satisfying $\ga(\tau)=M$ appeared in (\ref{6.9}).
For $\om\in\U$ and $M\in\Sp(2n)$, splitting numbers $S_N^{\pm}(\om)$ are constant
for all $N\in\Om^0(M)$. }
\\\\
{\bf Lemma 4.4.} (cf. \cite{Lon2}, Lemma 9.1.5 and List 9.1.12 of \cite{Lon5})
{\it For $M\in\Sp(2n)$ and $\om\in\U$, there hold
\bea
S_M^{\pm}(\om) &=& 0, \qquad {\it if}\;\;\om\not\in\sg(M).  \label{3.19}\\
(S_{N_1(1,a)}^+(1),S_{N_1(1,a)}^-(1))&=& \left\{\matrix{(1,1), &\quad {\rm if}\;\; a\ge 0, \cr
(0,0), &\quad {\rm if}\;\; a< 0. \cr}\right. \label{3.20}.\\
(S_{N_1(-1,a)}^+(-1),S_{N_1(-1,a)}^-(-1))&=& \left\{\matrix{(1,1), &\quad {\rm if}\;\; a\leq 0, \cr
(0,0), &\quad {\rm if}\;\; a> 0. \cr}\right. \label{3.21} \\
(S_{R(\theta)}^+(e^{\sqrt{-1}\theta}),S_{R(\theta)}^-(e^{\sqrt{-1}\theta}))&=&
(0,1) \;\ {\rm if}\;\ e^{\sqrt{-1}\theta}\in\U\setminus\R.\label{3.22}
\eea
If $e^{\sqrt{-1}\theta}\in\U\setminus\R$ and
$N_{2}(\omega,b)$ is nontrivial, then
\bea
(S_{N_2(\omega,b)}^+(e^{\sqrt{-1}\theta}),S_{N_2(\omega,b)}^-(e^{\sqrt{-1}\theta}))&=&
(1,1).\label{3.23}
\eea
If $e^{\sqrt{-1}\theta}\in\U\setminus\R$ and
$N_{2}(\omega,b)$ is trivial, then
\bea
(S_{N_2(\omega,b)}^+(e^{\sqrt{-1}\theta}),S_{N_2(\omega,b)}^-(e^{\sqrt{-1}\theta}))&=&
(0,0).\label{3.24}
\eea
For any $M_i\in\Sp(2n_i)$ with $i=0$ and $1$, there holds }
\be S^{\pm}_{M_0\dm M_1}(\om) = S^{\pm}_{M_0}(\om) + S^{\pm}_{M_1}(\om)
 \;\ {\rm and}\;\ S^{\pm}_{M}(\omega)=S^{\mp}_{M}(\bar{\omega}),   \;\ \forall\;\om\in\U. \lb{6.18}\ee
where $\bar{\omega}$ is the conjugate of $\omega$. Then we have the following
\\\\
{\bf Theorem 4.5.} (cf. \cite{Lon4} and Theorem 1.8.10 of \cite{Lon5}) {\it For
any $M\in\Sp(2n)$, there is a path $f:[0,1]\to\Om^0(M)$ such that $f(0)=M$ and
\be f(1) = M_1\dm\cdots\dm M_k,  \lb{6.19}\ee
where each $M_i$ is a basic normal form listed in (\ref{6.12})-(\ref{6.15})
for $1\leq i\leq k$.}

The following is the precise index iteration formulae for symplectic paths, which
is due to Y.Long (cf. Chaper 8 \cite{Lon5} or Theorem 2.1, 6.5 and 6.7 of \cite{LoZ1})
\\\\
{\bf Theorem 4.6.} For $n\in \N, \tau>0$, and any path $\gamma\in\P_{\tau}(2n)$,
set $M=\gamma(\tau)$. Extend $\gamma$ to the whole $[0,+\infty)$. Then for any
$m\in\N$,
\bea\nonumber
i(\gamma,m)&=&m(i(\gamma,1)+S^{+}_{M}(1)-C(M)) \\
&+&2\sum_{\theta\in(0,2\pi)}E(\frac{m\theta}{2\pi})S^{-}_{M}(e^{\sqrt{-1}\theta})
-(S_{M}^{+}(1)+C(M)). \label{3.23}
\eea
where $C(M)=\sum_{0<\theta<2\pi}S^{-}_{M}(e^{\sqrt{-1}\theta})$.
\\\\
{\bf Theorem 4.7.} Let $\ga\in\P_\tau(2n)$. Then there exists a path $f\in C([0,1],
\Omega^{0}(\ga(\tau)))$ such that $f(0)=\ga(\tau)$ and
\bea
f(1)&=&N_{1}(1,1)^{\diamond p_{-}}\diamond I_{2p_{0}} \diamond N_{1}(1,-1)^{\diamond p_{+}} \nonumber
\diamond N_{1}(-1,1)^{\diamond q_{-}}\diamond -I_{2q_{0}} \diamond N_{1}(-1,-1)^{\diamond q_{+}} \\ \nonumber
&\diamond& R(\theta_{1})\diamond\cdots R(\theta_{r}) \diamond N_{2}(\omega_{1},u_{1})\diamond
\cdots \diamond N_{2}(\omega_{r_{*}},u_{r_{*}}) \\
&\diamond& N_{2}(\lambda_{1},\upsilon_{1}) \diamond
\cdots \diamond N_{2}(\lambda_{r_{0}},\upsilon_{r_{0}})\diamond M_{k}.
\eea
where $N_{2}(\omega_{j},u_{j})$ are non-trivial form with some $\omega_{j}=e^{\sqrt{-1}\alpha_{j}}$,
$\alpha_{j}\in (0,\pi)\cup (\pi,2\pi)$ and $u_{j}=\left(
                                                        \begin{array}{cc}
                                                          u_{j1} & u_{j2} \\
                                                          u_{j3} & u_{j4} \\
                                                        \end{array}
                                                      \right)\in \R^{2\times2}
$,
$N_{2}(\omega_{j},u_{j})$
are trivial form with some  $\lambda_{j}=e^{\sqrt{-1}\beta_{j}}$,
 $\beta_{j}\in (0,\pi)\cup (\pi,2\pi)$ and $\nu_{j}=\left(
                                                        \begin{array}{cc}
                                                          \nu_{j1} & \nu_{j2} \\
                                                          \nu_{j3} & \nu_{j4} \\
                                                        \end{array}
                                                      \right)\in \R^{2\times2}
$, $M_{k}
=D(2)^{k}$ or $D(-2)\diamond D(2)^{\diamond(k-1)}$ ;
$p_{-}, p_{0}, p_{+}, q_{-}, q_{0}, q_{+}, r, r_{*}$ and $r_{0}$ are
non-negative integers; these integer and real number are uniquely determined by
$\ga(\tau)$. It holds that
\bea
n=p_{-}+p_{0}+p_{+}+q_{-}+q_{0}+q_{+}+r+2r_{*}+2r_{0}+k.\label{3.25}
\eea
We also have
$i(\ga,1)$ is odd if $f(1)=N_{1}(1,1),I_{2},N_{1}(-1,1),-I_{2},N_{1}(-1,-1)$
and $R(\theta)$; $i(\ga,1)$ is even if $f(1)=N_{1}(1,-1)$ and $N_{2}(\omega,b)$;
$i(\ga,1)$ can be any integer if $\sigma(f(1))\cap \U=\emptyset$. Then using the functions defined in (\ref{1.2}) , we have
\bea\nonumber
i(\ga,m)&=&m(i(\ga,1)+p_{-}+p_{0}-r)+2\sum_{j=1}^{r}E(\frac{m\theta_{j}}{2\pi})-r-p_{-}-p_{0} \\
&-&\frac{1+(-1)^{m}}{2}(q_{0}+q_{+})+2(\sum_{j=1}^{r_{*}}\varphi(\frac{m\alpha_{j}}{2\pi})-r_{*}) \label{3.26}
\eea
\bea \nonumber
\nu(\ga,m)&=&\nu(\ga,1)+\frac{1+(-1)^{m}}{2}(q_{-}+2q_{0}+q_{+})+2(r+r_{*}+r_{0}) \\
&-&2(\sum_{j=1}^{r}\varphi(\frac{m\theta_{j}}{2\pi})+\sum_{j=1}^{r_{*}}\varphi(\frac{m\alpha_{j}}{2\pi})+
\sum_{j=1}^{r_{0}}\varphi(\frac{m\beta_{j}}{2\pi}))  \label{3.27}
\eea
\bea
\hat{i}(\ga,1)=i(\ga,1)+p_{-}+p_{0}-r+\sum_{j=1}^{r}\frac{\theta_{j}}{\pi} \label{3.28}
\eea
\bea
S_{M}^{+}(1)=p_{-}+p_{0} \label{3.29}
\eea
\bea
C(M)=\sum_{0<\theta<2\pi}S_{M}^{-}(e^{\sqrt{-1}\theta})=q_{0}+q_{+}+r+2r_{*}.  \label{4.34}
\eea
\\\\

\medskip

\noindent {\bf Acknowledgements.}  The first author is  sincerely thanks Y.Long and C.Zhu for  the explanation of their methods and helpful discustion of this problem.

\bibliographystyle{abbrv}

\end{document}